\documentclass[11pt]{article}
\usepackage{amsmath,amssymb}
\usepackage{pstricks}
\usepackage{pst-node}
\usepackage[dvips]{graphicx}
\usepackage{pst-poly}
\usepackage{multido}

\newtheorem{result}{Theorem}

\newtheorem{problem}{Problem}
\newtheorem{define}{Definition}
\newtheorem{support}{Lemma}
\newtheorem{propo}{Proposition}

\newtheorem{note}{Remark}
\newtheorem{eg}{Example}

\newcommand{\qed}{%
\ifmmode 
\else \leavevmode\unskip\penalty9999 \hbox{}\nobreak\hfill \fi
\quad\hbox{\qedsymbol}}
\newcommand{\openbox}{\leavevmode \hbox to.77778em{%
\hfil\vrule
\vbox to.675em{\hrule width.6em\vfil\hrule}%
\vrule\hfil}}
\newcommand{\qedsymbol}{\openbox}

\newcommand{\showgrid}{}
\newcommand{\gridon}{\renewcommand{\showgrid}{\psset{subgriddiv=1,griddots=10,gridlabels=6pt}\psgrid}}
\gridon

\begin{document}

\begin{center} {\bf\LARGE  Triangulability  of Convex Graphs and Convex Skewness } \end{center}

\vskip5pt

\centerline{\large   Niran Abbas Ali$^a$, Gek L. Chia$^{b,e}$,  \  Hazim Michman  Trao$^{c}$ \ and \   Adem Kilicman$^d$ }

\begin{center}
\itshape\small  $^{a, c, d}\/$Department of Mathematics, \\ Universiti Putra Malaysia, 43400 Serdang, Malaysia,  \\
 \vspace{1mm}
 $^{b}\/$Department of Mathematical and Actuarial Sciences, \\  Universiti Tunku Abdul Rahman, Sungai Long Campus,    Malaysia \\
\vspace{1mm}
 $^e\/$Institute of Mathematical Sciences, University of Malaya, \\ 50603 Kuala Lumpur,  Malaysia  \\
\end{center}

\begin{abstract}
Motivated by  a result of \cite{aic:refer} which states that if  $F\/$ is a subgraph of a  convex complete graph $K_n\/$ and  $F\/$ contains no boundary edge of $K_n\/$  and $|E(F)| \leq n-3\/$, then $K_n - F\/$ admits a triangulation, we determine necessary and sufficient conditions on $F\/$ with $|E(F)| \leq n-1\/$ for which the conclusion remains true. For $|E(F)| \geq n\/$,  we investigate the possibility of packing $F\/$ in $K_n\/$ such that $K_n -F\/$ admits a triangulation for certain families of graphs $F\/$. These results are then applied to determine the convex skewness of the convex graphs of the form $K_n - F\/$.
\end{abstract}

\vspace{1mm}
 \section{Introduction and Preliminary}

By a {\em geometric graph\/} we mean a graph whose edges are straight line segments. By a {\em convex graph\/}, we mean a geometric graph whose vertices are in convex position. Let $G\/$ be a geometric graph with $n\/$ vertices having $k\/$ vertices in the convex hull.  By a {\em triangulation of $G\/$\/} we mean a maximal planar subdivision with the vertex set $V(G)\/$ of $G\/$. Hence a triangulation of $G\/$ has $2n-2-k\/$ triangles and $3n-3-k\/$ edges. It has been shown in \cite{llo:refer} that the problem of deciding whether a given geometric graph admits a triangulation of its vertex set is an $NP\/$-complete problem.

\vspace{1mm} For the case where $G\/$ is a convex geometric graph, it has been proved in \cite{aic:refer} that $G\/$  admits  a triangulation if $G\/$ is  obtained from the complete convex graph $K_n\/$ by deleting a set $F\/$ of at most $n-3\/$ edges and $F\/$ contains no edges from the boundary of the  convex hull of $G\/$ (see Theorem \ref{n-3}). It was noted that the result is best possible with respect to the number of edges in $G\/$. Motivated by this, we wish to obtain conditions on $F\/$ for which the result remains true when $F\/$ has more than $n-3\/$ edges.

\vspace{1mm} Throughout, if $F\/$ is  subgraph of $K_n\/$, then $K_n - F\/$ will denote the convex graph obtained from $K_n\/$ by deleting the set of  edges of $F\/$.

\begin{result}  (\cite{aic:refer} (Theorem 4.1))  \label{n-3}

Let $F\/$ be a subgraph of a  convex complete graph $K_n\/$. Suppose $F\/$ contains no boundary edge of $K_n\/$  and $|E(F)| \leq n-3\/$. Then $K_n - F\/$ admits a triangulation.
\end{result}

\vspace{1mm} For $F\/$ having at most $n-1\/$ edges, we obtain necessary and sufficient conditions on $F\/$ such that $K_n -F\/$ admits a triangulation. These conditions are given in Propositions \ref{fstars} and \ref{jstars}, and      Theorems  \ref{n-2} and  \ref{n-1}.

\vspace{1mm}
 The case where $F\/$ has at least $n\/$ edges seems to be a little complicated.  For this we turn to look for a configuration for $F\/$ that can be packed in $K_n\/$ so that $K_n - F\/$ admits a triangulation. If such a configuration exists for $F\/$, then we say that $F\/$ is {\em potentially triangulable in $K_n\/$}. Potentially triangulable graphs are considered in Section \ref{pt} where we (i) determine precisely the value of $n\/$ for which the $n\/$-cycle is potentially triangulable in $K_n\/$ (Theorem \ref{fnconnected}), and (ii) characterize all $2\/$-regular graphs  which are potentially triangulable in $K_n\/$ (Theorem \ref{disconnected}).

\vspace{1mm}
The potentially triangulable  problem is extended to the regular case in Section \ref{regular}. Here, while unable to  solve the general case, we consider the generalized Petersen graph.

\vspace{1mm}   We end the paper by  showing an application of these result to the problem of determining  the minimum number of edges to be deleted from a given convex graph so that the resulting graph is a convex plane graph.

\vspace{1mm}

Throughout, we shall adopt the following notations.  Unless otherwise stated, the vertices of a convex complete graph $K_n\/$ will be denoted by $v_0, v_1, \ldots, v_{n-1}\/$ in cyclic ordered. Also, unless otherwise specified, any operation on the subscript of $v_i\/$ is reduced modulo $n\/$.

 \vspace{1mm}
\begin{support} \label{ear}
Let $F\/$ be a subgraph of a convex complete graph $K_n\/$. Assume that $F\/$ has  at most $n-1\/$ edges and having no boundary edge of $K_n\/$. Then $K_n - F\/$ has a vertex $v_i\/$ such that $v_{i-1}v_{i+1}\/$ is not an edge of $F\/$.
\end{support}

\vspace{1mm}  \noindent
{\bf Proof:}   If  the lemma  is not true, then it implies that $v_{i-1}v_{i+1}, v_{i-2}v_i, v_iv_{i+2}  \in E(F)\/$, and recursively, this implies that $F\/$ is a spanning subgraph of $K_n\/$ with minimum vertex-degree at least $2\/$. But this implies that $|E(F)| \geq n\/$, a contradiction.   \qed

\section{Graphs with at most $n-2\/$ edges}  \label{n-1oeless}

We begin by describing a configuration for $F_n\/$ with $n\/$ vertices and $n-2\/$ edges such that $K_n - F\/$ admits no triangulation.

\vspace{2mm}
\begin{define}
 Let $F_n\/$ denote a subgraph of $K_n\/$ with no isolated vertices and having  $n-2\/$ edges. Let ${\cal F}_n(^*)\/$ denote a configuration of $ F_n\/$ on $K_n\/$  such that either
 $E({\cal F}_n(^*)) = \{ v_0v_i, v_1v_{n-1} \ | \ i =2, 3, \ldots, n-2\}\/$ or else ${\cal F}_n(^*)\/$ satisfies the following conditions:

 \vspace{1mm} (i) For  $ i =0, 1, 2, \ldots, k-1\/$,  $d(v_i) \geq 2\/$ and $v_{i-1}v_{i+1} \in E({\cal F}_n(^*))\/$ where $k\/$ is some  natural number with $2 \leq k \leq n-4\/$.

 \vspace{1mm} (ii) For $i =k, k+1, \ldots, n-1\/$, $d(v_i) =1\/$  and $v_i\/$ is adjacent to some non-pendant vertex in ${\cal F}_n(^*)\/$.

 \vspace{1mm} (iii) For any two pendant vertices $u, v\/$ such that $uv_i\/$ and $vv_j\/$ are crossed in ${\cal F}_n(^*)\/$, then $|i-j| =1\/$. 
\end{define}

\vspace{1mm}
 \begin{eg}  Let $F\/$ be a subgraph be a subgraph of $K_{10}\/$ with vertex set $\{v_0, v_1, \ldots, v_9\}\/$ and edge set $\{v_0v_2, v_0v_8, v_1v_3, v_1v_7, v_1v_9, v_3v_5, v_3v_6, v_2v_4\}\/$. Then $F\/$ is an example of a ${\cal F}_{10}(^*)\/$  configuration in $K_{10}\/$.  Here $k=4\/$.    \end{eg}

\vspace{1mm} Note that in the above definition of ${\cal F}_n(^*)\/$, if $w \in \{v_0, v_{k-1}\}\/$, then  $w\/$ is adjacent to some pendant vertex of ${\cal F}_n\/$.

\vspace{1mm}
\begin{propo} \label{fstars} Suppose $F_n\/$ is a subgraph of the convex complete graph $K_n\/$ with a configuration ${\cal F}_n(^*)\/$ as defined above.  Then
  $K_n - {\cal F}_n(^*)\/$ admits no triangulation for any natural number $ n \geq 4\/$.
\end{propo}

\vspace{1mm}  \noindent
{\bf Proof:} We prove this by induction on $n\/$. The result is clearly true if $4 \leq n \leq 6\/$. Hence we assume that $n \geq 7\/$ and that the result is true for all convex graphs $K_m - {\cal F}_m(^*)\/$ where $m < n\/$.

\vspace{1mm} By Lemma \ref{ear}, $K_n - {\cal F}_n(^*)\/$ contains a vertex $v_i\/$ such that $v_{i-1}v_{i+1} \not \in E({\cal F}_n(^*))\/$. By the definition of ${\cal F}_n(^*)\/$, we see that $d_{{\cal F}_n(^*)}(v_i)=1\/$.

\vspace{1mm} Let $F_{n-1} = {\cal F}_n (^*) - v_i\/$ and consider $K_{n-1} - F_{n-1}\/$. We shall show that $F_{n-1}\/$ is a subgraph of $K_{n-1}\/$ with the configuration ${\cal F}_{n-1}(^*)\/$.

\vspace{1mm} Let $v_j\/$ be the neighbor of $v_i\/$ in $F_n\/$.

\vspace{2mm} {\em Case (1):  $d_{F_{n-1}}(v_{j}) \geq 2 $.   }

\vspace{1mm}
    Clearly $F_{n-1} = {\cal F}_{n-1}(^*)$  if $ i \geq j+3\/$ (otherwise $F_{n} \neq {\cal F}_n(^*)$). So assume that $i=j+2$.

 \vspace{1mm} If  $d_{F_{n}}(v_{j+1}) =1 $, then $v_j = v_{k-1}\/$ and  we have  $F_{n-1} = {\cal F}_{n-1}(^*)$ by definition.

 \vspace{1mm} If $d_{F_{n}}(v_{j+1}) \geq 2  $, then $v_{j+1} =v_{k-1}\/$ (since $v_{j+2} = v_i\/$ is a pendant vertex). Hence $v_{j+1}\/$ is adjacent to some pendant vertex $v_r\/$ of $F_n\/$ where $r \in \{j+4, \ldots, n-1\}\/$.

  \vspace{1mm} Suppose   $v_jv_{j+3} \not \in E(F_{n-1})\/$.  Then
   $v_{j+3}\/$ is adjacent to some vertex $v_s\/$ where $s \in \{0,1, \ldots, j-1\}\/$ with $d_{F_n}(v_s) \geq 2\/$. But this contradicts condition (iii) (in the definition of ${\cal F}_n(^*)\/$) since $v_{j+3}v_s\/$ crosses $v_rv_{j+1}\/$ and  $s \neq j\/$.

   \vspace{1mm} Hence $v_jv_{j+3} \in E(F_{n-1})\/$ and it follows that $F_{n-1} = {\cal F}_{n-1}(^*)\/$ in this case.

\vspace{2mm} {\em Case (2):  $d_{F_{n-1}}(v_j) =1 $.   }

\vspace{1mm} In this case, $d_{F_n}(v_j) =2\/$. We assert that $j \in \{0, k-1\}\/$.

\vspace{1mm} Suppose $1 \leq j \leq k-2\/$.  Since $v_{j-1}, v_{j+1}\/$ are non-pendant vertices of ${\cal F}_n(^*)\/$, it follows that  $v_{j-2}, v_{j+2}\/$ are the only neighbors of $v_j\/$ in ${\cal F}_n(^*)\/$. Assume without loss of generality that $v_i = v_{j+2}\/$. Hence $v_{j+1} = v_{k-1}\/$ is adjacent to a pendant vertex $v_r\/$ where $r \in \{j+4, \ldots, n-1\}\/$. Now the pendant vertex $v_{j+3}\/$ is adjacent to some non-pendant vertex $v_s\/$ where $s \in \{0, 1, \ldots,  j-1\}\/$. But then $v_{j+3}v_s, v_rv_{j+1}\/$ are crossing each other and $s \neq j\/$, a contradiction to condition (iii) in the definition of ${\cal F}_n(^*)\/$. This proves the assertion.

\vspace{1mm} Assume without loss of generality that $j =k-1\/$.  Then clearly $F_{n-1} = {\cal F}_{n-1}(^*)\/$.

\vspace{1mm}
By induction $K_{n-1} - {\cal F}_{n-1}(^*)\/$ admits no triangulation.

\vspace{1mm}
Assume on the contrary that $K_n - {\cal F}_n(^*)$ admits a triangulation $T_n\/$.

\vspace{1mm} (i) Suppose $T_n\/$ contains $v_{i-1}v_{i+1}\/$. Then $v_i\/$ is a degree-$2\/$ vertex in  $T_n\/$. Hence $T_n -v_i\/$ is a triangulation for $K_{n-1} - {\cal F}_{n-1}(^*)\/$, a contradiction.

\vspace{1mm} (ii) Suppose $T_n\/$ does not contain $v_{i-1}v_{i+1}\/$. Then $v_iv_t\/$ is a diagonal in   $T_n\/$ and $v_iv_t\/$ divides $K_n - F_n(^*)$ into $K_{m_1} - F_{m_1}$  and    $K_{m_2} - F_{m_2}$ where $m_1 + m_2 -2 = n\/$.  Without loss of generality assume that $v_iv_j\/$ is an edge in $F_{m_1}\/$. Then it is readily checked that $F_{m_1} = {\cal F}_{m_1}(^*)\/$. By induction $K_{m-1} - {\cal F}_{m_1}(^*)\/$ admits no triangulation, a contradiction.

\vspace{1mm} This completes the proof.  \qed

\vspace{2mm} The following result together with Proposition \ref{fstars} characterize all graphs $F_n\/$ with at most $n-2\/$ edges such that $K_n - F_n\/$ admits a triangulation.

\vspace{1mm}
\begin{result}    \label{n-2}
Suppose $n \geq 6\/$ is a natural number and $F_n\/$ is a subgraph of the  convex complete graph $K_n\/$ such that  $|E(F_n)| \leq n-2\/$ and $F_n\/$ contains no boundary edges  of $K_n\/$. Then $K_n - F_n\/$ admits a triangulation unless $F_n = {\cal F}_n(^*) \/$.
\end{result}

\vspace{3mm}  \noindent
{\bf Proof:} In view of Proposition \ref{fstars}, we assume that $F_n \neq {\cal F}_n(^*) \/$.

\vspace{1mm}
We prove the result by induction on $n\/$.   The result is clearly true if $n =6\/$. Assume that $n \geq 7\/$ and the result is true for all convex graphs $K_m - F_m\/$ where $m < n\/$.

\vspace{1mm} By Lemma \ref{ear}, $K_n - F_n\/$ contains a vertex $v_i\/$ such that $v_{i-1}v_{i+1} \not \in E(F_n)\/$.  \hfill $(*)$

\vspace{1mm} Let $v_i\/$ which satisfies the condition in $(*)\/$ be such that $d_{F_n}(v_i) = k\/$ is maximal.

\vspace{1mm}
Delete the  vertex $v_i\/$ from $K_n - F_n\/$. The resulting graph  is a convex graph of the form $K_{n-1} - F_{n-1}\/$ where $F_{n-1} = F_n -v_i\/$ is a subgraph of $K_{n-1}\/$.

\vspace{1mm}
{\em Case (1)}:  $k \geq 2\/$.

\vspace{1mm}
 Then $|E(F_{n-1})| \leq n-4 $.    Clearly  $F_{n-1}\/$ does not admit the configuration ${\cal F}_{n-1}(^*)$. By induction $K_{n-1} - F_{n-1}\/$ admits a triangulation $T_{n-1}\/$. As such $T_{n-1} \cup \{v_{i-1}v_i, v_i v_{i+1}\}\/$ is a triangulation for $K_n - F_n\/$.

\vspace{1mm}
{\em Case (2)}:  $k =1\/$.

\vspace{1mm}
Let $v_j\/$ be the neighbor of $v_i\/$ in $F_n\/$.

\vspace{1mm}
If $d_{F_n}(v_j) =1\/$, then $v_j\/$ is an isolated vertex in $F_{n-1}\/$ which implies that $K_{n-1} - F_{n-1}\/$ admits a triangulation $T_{n-1}\/$. Again $T_{n-1} \cup \{v_{i-1}v_i, v_i v_{i+1}\}\/$ is a triangulation for $K_n - F_n\/$.

\vspace{1mm}
Hence assume that $d_{F_n}(v_j) \geq 2\/$.

\vspace{1mm} If $F_{n-1} \neq {\cal F}_{n-1}(^*)\/$, then $K_{n-1} - F_{n-1}\/$ admits a triangulation $T_{n-1}\/$ by induction. Again $T_{n-1} \cup \{v_{i-1}v_i, v_i v_{i+1}\}\/$ is a triangulation for $K_n - F_n\/$.

\vspace{1mm}
Hence we assume that  $F_{n-1} = {\cal F}_{n-1}(^*)\/$.

\vspace{1mm}
In view of the maximality on $d_{F_n}(v_i) =k\/$, we see that $v_{i-1}, v_{i+1}\/$ are both  pendant vertices in $F_n\/$ (otherwise $F_n  = {\cal F}_n(^*)\/$). Let $v_s\/$ and $v_t\/$ be the neighbors of $v_{i-1}\/$ and $v_{i+1}\/$ in $F_n\/$ respectively.

\vspace{1mm}
If $v_iv_j\/$ crosses neither $v_{i-1}v_s\/$ nor $v_{i+1}v_t\/$, then $F_n  = {\cal F}_n(^*)\/$,  a contradiction.

\vspace{1mm} Hence we assume without loss of generality that $v_iv_j\/$ crosses $v_{i-1}v_s\/$ and that $s < j < i\/$.

\vspace{1mm} If $|s-j| =1\/$, then we have $F_n = {\cal F}_n(^*)\/$, a contradiction.

\vspace{1mm} If $|s-j| > 1\/$,
  a triangulation for $K_n - F_n\/$ is given by the set of diagonals $v_{i-1}v_{j-1}\/$,  $v_{i-1}v_l\/$ where $ l \in \{   j , j+1, \ldots, i-3\} =A  \/$ and $v_iv_m\/$ for all $m \in \{0, 1, 2, \ldots, n-1\} - (A \cup \{i-2, i-1, i, i+1\})\/$.

\vspace{1mm} This completes the proof.   \qed

\section{Graphs with only $n-1\/$ edges}  \label{n-1}

\vspace{2mm} We  begin by describing  a  configuration for $F_n\/$ with $n\/$ vertices and $n-1\/$ edges such that $K_n - F\/$ admits no triangulation.

\vspace{2mm}
\begin{define}  \label{j2}
 Let $F_n\/$ denote a subgraph of $K_n\/$ with no isolated vertices and having  $n-1\/$ edges where $n \geq 5\/$. Let ${\cal J}_n(^*)\/$ denote a configuration of $ F_n\/$ on $K_n\/$  such that (i) whenever $v_{i-1} v_{i+1}\not \in E({\cal J}_n(^*))\/$, then $d_{F_n}(v_i) \leq 2\/$, and (ii) in the case that $d_{F_n}(v_i) =2\/$, then ${\cal J}_n(^*) -v_i = {\cal F}_{n-1}(^*)\/$ in $K_n - v_i\/$.
\end{define}

\vspace{2mm}
Suppose $k \in \{1, 2\}\/$. Let  $J_{n,k}\/$   denote the set of all degree-$k\/$ vertices in ${\cal J}_n(^*)\/$ satisfying the conditions (i) and (ii) in Definition \ref{j2}. Note that if $v \in J_{n,2}\/$, then no neighbor of $v\/$ is a pendant vertex in ${\cal J}_n(^*)\/$.

\begin{define}
The configuration ${\cal J}_n(^*)\/$  is said to be of {\em Type-1\/} if $J_{n,2} = \emptyset\/$ and there is an edge $e\/$ in  ${\cal J}_n(^*)\/$  such that ${\cal J}_n(^*)-e\/$  is ${\cal F}_n(^*)\/$.  The configuration ${\cal J}_n(^*)\/$  is said to be of {\em Type-2\/} if (i) $J_{n,2} \neq \emptyset\/$ and (ii) whenever $v_i \in J_{n,2}\/$ and $v_t\/$ is a pendant vertex in ${\cal J}_n(^*)\/$ such that $v_tv_s\/$ crosses $v_iv_{j_1}\/$ and $v_iv_{j_2}\/$,  then either $|j_1 -s| =1\/$ or $|j_2 -s |=1\/$.
\end{define}

\vspace{1mm}
 \begin{eg}
 (i) Let   $F\/$ be a subgraph of $K_{11}\/$ with vertex set $\{ v_0, v_1, \ldots, v_{10}\}\/$ and edge set $\{v_0v_2, v_0v_3, v_0v_7, v_0v_9, v_1v_3, v_1v_8, v_1v_{10}, v_2v_4, v_2v_6, v_3v_5\}\/$. Then $F\/$ is an example of a  ${\cal J}_{11}(^*)\/$ configuration in $K_{11}\/$. Here  $J_{11,1} = \{v_5, v_6, \ldots, v_9\}\/$ and $J_{11,2} = \emptyset\/$. Hence $F\/$ is of Type-$1\/$.  Note that ${\cal J}_{11}(^*) - \{v_0v_3\}  \cong {\cal F}_{11}(^*)\/$.

 \vspace{1mm}
 (ii)  Let $F\/$  be a subgraph of $K_{11}\/$ with vertex set $\{ v_0, v_1, \ldots, v_{10}\}\/$ and edge set $\{v_0v_2, v_0v_7, \linebreak  v_0v_9, v_1v_3, v_1v_8, v_1v_{10}, v_2v_4, v_2v_6, v_3v_5,  v_5v_9\}\/$. Then $F\/$ is an example of a ${\cal J}_{11}(^*)\/$ configuration in $K_{11}\/$. Here  $J_{11,1} = \{ v_6, v_7, v_8\}\/$ and $J_{11,2} = \{v_5, v_9\}\/$. Hence $F\/$ of  Type-$2\/$.          \end{eg}

\vspace{1mm}
We shall define three  ${\cal J}_6(^*)\/$ configurations $J_1, J_2, J_3\/$ (on $K_6\/$) each  of Type-$2$ with $J_{6,1} = \emptyset$. Let $J_1 =  v_1v_4v_2v_5v_1 \cup \{v_0v_3\}\/$, $J_2 = v_1v_3v_5v_1 \cup \{v_0v_2, v_1v_4\}\/$ and $J_3 = v_1v_3v_5v_1 \cup v_2v_0v_4\/$.

\begin{support} \label{n=6}
 Suppose $F_6$ is a Type-$2$ ${\cal J}_6(^*)\/$. Then  $J_{6,1} = \emptyset\/$ if and only if $F_6 = J_i\/$ for some $i \in \{1, 2, 3 \}$. Moreover  $K_6 -J_i$ admits no triangulation for any $i =1, 2, 3\/$.
\end{support}

\vspace{1mm}  \noindent {\bf Proof:} The sufficiency is clear. We shall prove the necessity.

\vspace{1mm} Note that, for any  $v \in J_{6,2}$, ${\cal J}_6(^*) - v = {\cal F}_5(^*)\/$. Note also that there is  only one  configuration of  ${\cal F}_5(^*)\/$  on $K_5$ which is given by  $ E({\cal F}_5(^*))= \{x_0x_2, x_1x_3, x_1x_4\}\/$.

\vspace{3mm} Now, there are only three possible ways to place $v\/$ on the boundary edge of $K_5\/$ (in order to recover ${\cal J}_6(^*)\/$ from ${\cal F}_5(^*)\/$). Without loss of generality, assume that $v\/$ lies on the edge $x_ix_{i+1}\/$ for some $i \in \{1, 2, 3\}\/$. In each case, we join $v\/$ to two  vertices of ${\cal F}_5(^*)\/$ subject to the condition that $J_{6,1} = \emptyset\/$. We then arrive at the conclusion that $F_6 = J_i\/$, $i \in \{1, 2, 3\}\/$.

\vspace{1mm} Suppose on the contrary that $K_6- J_i\/$ admits a triangulation $T_i\/$. For any $v_j \in J_{6,2}\/$, the edge $v_{j-1}v_{j+1}\/$ is not a diagonal edge in $T_i\/$ (because $J_i - v_j = {\cal F}_5(^*)\/$).  Hence $v_jv_t\/$, for some vertex $v_t\/$ (not adjacent to $v_j\/$) in $J_i\/$,  is a diagonal edge in  $T_i\/$. But then it is routine to check that $v_jv_t\/$ cannot be combined with another two non-adjacent edges in $J_i\/$ to form the set of diagonal edges of $T_i\/$, a contradiction.    \qed

\begin{support} \label{j1emety}
  Suppose $F_n$ is a Type-$2$ ${\cal J}_n(^*)\/$ configuration and $n\geq 7$,  If $J_{n,1} = \emptyset$, then for each $v_i \in J_{n,2}$, $v_{i-2}v_i, v_iv_{i+2} \in E({\cal J}_n(^*)\/)$.
\end{support}

\vspace{1mm}  \noindent {\bf Proof:}
Since  $v_i \in J_{n,2}$, we have  ${\cal J}_n(^*) - v_i = {\cal F}_{n-1}(^*)\/$. Let $u_0,  u_1, \ldots, u_{n-2}\/$ denote the vertices of  ${\cal F}_{n-1}(^*) \/$   arranged in cyclic order.

\vspace{1mm}
Now,   $J_{n,1}= \emptyset$ implies that   ${\cal F}_{n-1}(^*)\/$ has  at least two pendant vertices and at most three pendant vertices $u_{j}$  such that $u_{j -1}u_{j +1} \not \in   E({\cal F}_{n-1}(^*))$. Moreover such pendant vertices must be in consecutive order (in ${\cal F}_{n-1}(^*)\/$).

 \vspace{1mm} Now, insert the vertex $v_i\/$ into  ${\cal F}_{n-1}(^*)\/$  to recover ${\cal J}_n(^*)\/$. Since $n \geq 7\/$ and  $J_{n,1} = \emptyset$,  $v_i\/$ must be inserted in between two pendant vertices of ${\cal F}_{n-1}(^*)\/$.   

 \vspace{1mm} If there are only two such pendant vertices $u_j \/$ and $u_{j+1}\/$, then either $u_j = v_{i-1}, \ u_{j+1} = v_{i+1}\/$
 or else  $u_j = v_{i+1}, \ u_{j-1} = v_{i-1}\/$. Either case implies that   $v_{i-2}v_i, v_iv_{i+2} \in E({\cal J}_n(^*)\/)$.

 \vspace{1mm} If there are three such pendant vertices $u_{j-1}\/$, $u_j \/$ and $u_{j+1}\/$, then again either  $u_j = v_{i-1}, \ u_{j+1} = v_{i+1}\/$
 or else  $u_j = v_{i+1}, \ u_{j-1} = v_{i-1}\/$ and we have the same conclusion as before.

 \vspace{1mm} This proves the lemma.  \qed.

\begin{propo} \label{jstars} Suppose $F_n\/$ is a subgraph of the convex complete graph $K_n\/$ with a Type-$1\/$ or Type-$2\/$    configuration ${\cal J}_n(^*)\/$.  Then
  $K_n - {\cal J}_n(^*)\/$ admits no triangulation for any natural number $ n \geq 5\/$.
\end{propo}

\vspace{1mm}  \noindent {\bf Proof:}  We prove this by induction on $n\/$. The result is clearly true if $ n =5\/$. Hence we assume that $n \geq 6\/$ and that the result is true for all convex graphs $K_m - F_m\/$ where $m < n\/$.

\vspace{1mm} By Lemma \ref{ear}, $ K_n - {\cal J}_n(^*) \/$ contains a vertex $v_j\/$ such that $v_{j-1}v_{j+1} \not \in E({\cal  J}_n(^*))\/$. By the definition of ${\cal  J}_n(^*) \/$, we see that $d_{{\cal  J}_n(^*)}(v_j)\leq 2\/$.

\vspace{2mm} {\em Case (1): ${\cal  J}_n(^*)\/$ is of Type-$1\/$.   }

\vspace{1mm}
In this case $d_{{\cal  J}_n(^*)}(v_j)=1\/$, and there  is an edge  $e \in {\cal  J}_n(^*)$ satisfying ${\cal  J}_{n}(^*)-e \cong {\cal F}_{n}(^*)\/$. Since $K_n - {\cal F}_n(^*)$ admits no  triangulation, it follows that  $K_n - {\cal  J}_n(^*)$ admits no  triangulation

\vspace{2mm} {\em Case (2):   ${\cal  J}_n(^*)\/$ is of Type-$2\/$.   }

\vspace{1mm}
(a) Suppose $J_{n,1} \neq \emptyset\/$.

\vspace{1mm} In this case, let $v_j\/$ which satisfies the condition $v_{j-1}v_{j+1} \not \in E({\cal  J}_n(^*))\/$  be chosen such that   $d_{{\cal  J}_n(^*)}(v_j)=1\/$.

\vspace{1mm} Then ${\cal  J}_n(^*)  - v_j = {\cal J}_{n-1}(^*)\/$  in $K_{n-1}\/$. By induction, $K_{n-1} - {\cal  J}_{n-1}(^*)  \/$ admits no triangulation.

\vspace{1mm}
Assume on the contrary that $K_n - {\cal  J}_n(^*)$ admits a triangulation $T_n\/$.

\vspace{1mm}  Suppose $T_n\/$ contains $v_{j-1}v_{j+1}\/$. Then  $T_n -v_j\/$ is a triangulation for $K_{n-1} - {\cal  J}_{n-1}(^*)  \/$, a contradiction.

\vspace{1mm} Hence assume that $T_n\/$ does not contain $v_{j-1}v_{j+1}\/$. Then $T_n\/$ contains a diagonal $v_jv_t\/$. Let $v_{j_1}$  be the neighbor of $v_j$ in ${\cal  J}_n(^*)$.

 \vspace{1mm}
Clearly $v_jv_t\/$ divides $K_n - {\cal  J}_n(^*)$ into two convex graphs  $K_{m_1} - F_{m_1}$  and    $K_{m_2} - F_{m_2}$ where $m_1 + m_2  = n+2\/$. Without loss of generality assume that $v_jv_{j_1}\/$ is an edge in $F_{m_1}\/$. Note that $K_{m_2} - F_{m_2}\/$ admits a triangulation (since $v_j\/$ is an isolated vertex in  $ F_{m_2}\/$).

\vspace{1mm}
 If $J_{n,2}\cap V(F_{m_1}) = \emptyset\/$,  then $F_{m_1} = {\cal F }_{m_1}(^*)$ which admits no triangulation by Proposition \ref{fstars}. Hence we assume that $J_{n,2} \cap V(F_{m_1})\/$ contains a vertex $v_i\/$.

\vspace{1mm}
(i) If  both neighbors of $v_i\/$ are in $V(F_{m_1})\/$, then  $F_{m_1}\/$ is ${\cal  J}_{m_1}(^*)$. By induction, $K_{m_1} - F_{m_1}\/$ admits no triangulation and this implies that $K_n - {\cal  J}_n(^*)$ admits no triangulation, a contradiction.

\vspace{1mm} (ii) If only one of the neighbors of $v_i\/$ is in $V(F_{m_1})\/$, then $F_{m_1}\/$ is ${\cal  F}_{m_1}(^*)$ and again we have a contradiction.

\vspace{1mm}
(iii) If  both neighbors $v_{i_1}, v_{i_2}\/$ of $v_i\/$ are in $V(F_{m_2})\/$, then the edge $v_jv_s\/$ (incident to the pendant vertex $v_j\/$) crosses both the edges $v_iv_{i_1}, v_iv_{i_2}\/$  with $|s-i_1|>1\/$ and $|s - i_2| > 1\/$, a contradiction to the definition of Type-$2\/$  ${\cal  J}_n(^*)$ configuration.

\vspace{1mm}
(b) Suppose  $ J_{n,1}= \emptyset\/$.

 \vspace{1mm} Then $v_j \in J_{n,2}\/$ and ${\cal  J}_n(^*) - v_j = {\cal F}_{n-1}(^*)\/$  in $K_{n-1}\/$. By Proposition \ref{fstars}, $K_{n-1} - {\cal  F}_{n-1}(^*)  \/$ admits no triangulation.

\vspace{1mm}   When $n=6$,   $K_6 - {\cal  J}_6(^*)\/$ admits  no triangulation  by Lemma \ref{n=6}.

\vspace{1mm} When $n\geq7 $, assume on the contrary that $K_n - {\cal  J}_n(^*)$ admits a triangulation $T_n\/$.

\vspace{1mm}  If $T_n\/$ contains $v_{j-1}v_{j+1}\/$, then  $T_n -v_j\/$ is a triangulation for $K_{n-1} - {\cal F}_{n-1}(^*)  \/$, a contradiction.

\vspace{1mm} Hence assume that $T_n\/$ does not contain $v_{j-1}v_{j+1}\/$. Then $T_n\/$ contains a diagonal $v_jv_t\/$.

\vspace{1mm} By Lemma \ref{j1emety}, $v_{j-2}v_j,  \ v_jv_{j+2} \in  E({\cal J}_n(^*)\/)$ (since $n \geq 7\/$).

\vspace{1mm} Clearly $v_jv_t\/$ divides $K_n - {\cal  J}_n(^*)$ into two convex graphs  $K_{m_1} - F_{m_1}$  and    $K_{m_2} - F_{m_2}$ where $m_1 + m_2 -2 = n\/$ with $F_{m_1}\/$  and $F_{m_2}\/$ containing  $v_jv_{j-2}$ and $v_j v_{j+2}\/$ respectively.

 \vspace{1mm}
  Note that $v_t \not \in J_{n,2}\/$. This is because   the vertex  $v_j$  in ${\cal  J}_n(^*) - v_t\/$ is a non-pendant vertex and    $v_{j-1}v_{j+1}\/$ is not an edge of ${\cal J}_{n}(^*) - v_t$ implying that ${\cal  J}_n(^*) - v_t \neq {\cal  F}_{n-1}(^*)$.

\vspace{1mm} If  $V(F_{m_1}) \cap J_{n,2} = \emptyset\/$,    then clearly $F_{m_1}\/$ is ${\cal  F}_{m_1}(^*)$. By Proposition \ref{fstars}, $K_{m_1} - {\cal  F}_{m_1}(^*)\/$  admits no triangulation and this implies that $K_n - {\cal  J}_n(^*)$ admits no triangulation, a contradiction.

\vspace{1mm} Hence assume that there exists $v_i \in V(F_{m_1}) \cap J_{n,2}\/$. Clearly, $v_i$ is adjacent to $v_j$ in ${\cal  J}_n(^*)$ (otherwise, the vertex  $v_j$  in ${\cal  J}_n(^*) - v_i$ is a non-pendant vertex with  $v_{j-1}v_{j+1} \notin ({\cal J}_{n}(^*) - v_i)$, and this implies that ${\cal  J}_n(^*) - v_i \neq {\cal  F}_{n-1}(^*)$, a contradiction).   Thus  $i = j-2$.

\vspace{1mm} We assert that $(F_{m_2} - v_j) \cap J_{n,2} = \emptyset\/$.   To see this,    suppose  $v_r \in (F_{m_2} - v_j) \cap J_{n,2} \/$. By the same argument in the preceding paragraph (on $v_i\/$), it follows that $v_r\/$ is adjacent to $v_j\/$ and hence $v_r = v_{j+2}\/$. Moreover  $v_r\/$ is not adjacent to $ v_{i}\/$ (since $n \geq 7\/$). As such,  ${\cal  J}_n(^*) - v_{r}$  contains  $v_i$  with $ v_{i-1}v_{i+1} \notin {\cal J}_{n}(^*) - v_{r}$ (and $v_i\/$ is non-pendant)  and this implies that $ {\cal J}_{n}(^*) - v_i \neq  {\cal F}_{n-1}(^*)$, a contradiction.

\vspace{1mm} It follows from the assertion that  $F_{m_2} = {\cal  F}_{m_2}(^*)$. By Proposition \ref{fstars},  $K_{m_2} - {\cal  F}_{m_2}(^*)\/$ admit no triangulation and this implies that $K_n - {\cal  J}_n(^*)$ admits no triangulation, a contradiction.

\vspace{1mm} This completes the proof.   \qed

\vspace{1mm}
\begin{result}    \label{n-1}
Suppose $n \geq 5\/$ is a natural number and $F_n\/$ is a subgraph of the  convex complete graph $K_n\/$ such that  $|E(F_n)| = n-1\/$ and $F_n\/$ contains no boundary edges  of $K_n\/$. Then $K_n - F_n\/$ admits a triangulation unless $F_n \/$ is Type-1 or Type-2 ${\cal J}_n(^*)\/$.
\end{result}

\vspace{3mm}  \noindent
{\bf Proof:} We prove the result by induction on $n\/$.

\vspace{1mm}
 The result is clearly true if $n = 5, 6\/$. Assume that $n \geq 7\/$ and that  the result is true for all convex graphs $K_m - F_m\/$ where $m < n\/$.

\vspace{1mm} By Lemma \ref{ear}, $K_n - F_n \/$ contains a vertex $v_i\/$ such that $v_{i-1}v_{i+1} \not \in E(F_n)\/$.  \hfill $(*)$

\vspace{1mm} Let $v_i\/$ which satisfies the condition in $(*)\/$ be such that $d_{F_n}(v_i) = k\/$ is maximal.

\vspace{1mm}
Delete the  vertex $v_i\/$ from $K_n - F_n\/$. The resulting graph  is a convex graph of the form $K_{n-1} - F_{n-1}\/$ where $F_{n-1} = F_n -v_i\/$ is a subgraph of $K_{n-1}\/$.

\vspace{2mm}
{\em Case (1)}:  $k \geq 3\/$.

\vspace{1mm}
 Then $|E(F_{n-1})| \leq n-4 $.    By Theorem  \ref{n-2},  $K_{n-1} - F_{n-1}\/$ admits a triangulation $T_{n-1}\/$. As such $T_{n-1} \cup \{v_{i-1}v_i, v_i v_{i+1}\}\/$ is a triangulation for $K_n - F_n\/$.

\vspace{2mm}
{\em Case (2)}:  $k =2\/$.

\vspace{1mm} If $F_{n-1} \neq {\cal F}_{n-1}(^*)\/$, then $K_{n-1} - F_{n-1}\/$ admits a triangulation $T_{n-1}\/$ by Theorem 2. Again $T_{n-1} \cup \{v_{i-1}v_i, v_i v_{i+1}\}\/$ is a triangulation for $K_n - F_n\/$.

\vspace{1mm}Hence we assume that  $F_{n-1} = {\cal F}_{n-1}(^*)\/$.

\vspace{1mm}
Let $v_{j_1}\/$ and $v_{j_2}\/$ be the neighbors of $v_i\/$ in $F_n\/$.

\vspace{1mm}
Now,  there is an edge $v_tv_s\/$ incident to a  pendant vertex $v_t\/$ such that  $v_tv_s\/$ crosses both the edges $v_iv_{j_1}, \  v_iv_{j_1}\/$  with  $|j_1 -s| >1\/$ and  $|j_2 -s |> 1\/$ (otherwise $F_n\/$ is Type-2 ${\cal J}_n(^*)\/$).

\vspace{1mm} We can assume without loss of generality that   $s < j_1 < j_2 < i\/$.

\vspace{1mm} Let $G_1\/$ and $G_2\/$ be subgraphs of $K_n - F_n$ induced by the vertices $\{v_{j_1 -1}, v_{j_1}, \ldots,  v_{i-1}\}$ and $\{v_{i -1}, v_{i}, \ldots,   v_{n-1}, v_0,  \ldots, v_{j_1 -1}\}$  respectively (with $v_t \/$ contained in $G_1\/$).  Clearly,  $d_{G_1}(v_{t})=0$ and $d_{G_2}(v_i)=0$ which imply that  $G_i$ admits a triangulation $T_i$, $i=1,2$. As such,  $T= T_1 \cup T_2$ is a triangulation for  $K_n - F_n$.

\vspace{2mm}
{\em Case (3)}:  $k =1\/$.

\vspace{1mm}
Let $v_j\/$ be the neighbor of $v_i\/$ in $F_n\/$.

\vspace{1mm}
Then  $d_{F_n}(v_j) \geq 2\/$ otherwise  $K_{n-1} - F_{n-1}\/$ admits a triangulation $T_{n-1}\/$ (because $v_j\/$ is an isolated vertex in  $K_{n-1} - F_{n-1}\/$)  which means that $T_{n-1} \cup \{v_{i-1}v_i, v_i v_{i+1}\}\/$ is a triangulation for $K_n - F_n\/$.

\vspace{1mm}
By the  maximality of $k\/$, we have  $v_{j-1}v_{j+1} \in E(F_n)$.

\vspace{1mm} Now, $|E(F_{n-1})|= n-2$. If $F_{n-1}\/$ is not the configuration ${\cal J}_{n-1}(^*)\/$, then $K_{n-1} - F_{n-1}\/$ admits a triangulation $T_{n-1}\/$ by induction. Again $T_{n-1} \cup \{v_{i-1}v_i, v_i v_{i+1}\}\/$ is a triangulation for $K_n - F_n\/$.

\vspace{1mm}
Hence we assume that  $F_{n-1} \cong {\cal J}_{n-1}(^*)\/$.

\vspace{1mm}
(i)  Suppose  ${\cal J}_{n-1}(^*)\/$ is of Type-2.

 \vspace{1mm} Then there is a unique degree-$2\/$ vertex $v_{s} $ in  ${\cal J}_{n-1}(^*)\/$  such that $ {\cal J}_{n-1}(^*) - v_s =  {\cal F}_{n-2}(^*) \/$ (in $K_{n-2} = K_n - \{v_i, v_s\}\/$). It is easy to see that either $v_s = v_{i+1}\/$ or $v_s = v_{i-1}\/$.  Assume without loss of generality that $v_s= v_{i+1}$. Hence  $v_{i-1}v_{s+1}$ is not an edge  in $ {\cal J}_{n-1}(^*)\/$ (since  $v_s\in J_{n-1,2}$).

\vspace{1mm}
In view of maximality of $k$, $v_{s-1}v_{s+1} \in E(F_{n})\/$.  Since $v_i= v_{s-1}\/$ is a pendant vertex, it follows that  $v_j = v_{s+1}$.

\vspace{1mm} Clearly, $v_{i-1}\/$ is not adjacent to $v_{j}= v_{s+1}\/$. Hence   $v_{i-1} v_j $  is a boundary edge in $K_{n-2} -  {\cal F}_{n-2}(^*)$.

\vspace{1mm}
Let $v_{s_1}, v_{s_2}\/$ be the neighbors of $v_s\/$ in ${\cal J}_{n-1}(^*)\/$.

\vspace{1mm}
Since  $v_{j-1}v_{j+1} = v_sv_{s+2}\/$ is an edge in $F_n\/$, we may assume that    $v_{s+2}= v_{s_1}\/$, and hence $s_2 \geq s+3\/$.

\vspace{1mm}  By the maximality of $k\/$ we see that $v_{i-1}\/$ is a pendant vertex in   $F_n$ (since $v_i\/$ is a pendant vertex and $v_i\/$ is adjacent only to $v_j = v_{s+1}\/$).
 But this means that $v_{i-1}$ is also a  pendant vertex in ${\cal F}_{n-2}(^*)$.

\vspace{1mm}   Note that   $v_{s_1}\/$ is a non-pendant vertex in $F_n\/$ (otherwise $ {\cal J}_{n-1}(^*) - v_s \/$ contains an isolated vertex $v_{s_1}\/$, a contradiction).

\vspace{1mm} It is easy to see that ${\cal F}_{n-2}(^*) \cup \{v_{s}v_{s_1}\} \cup \{v_iv_j\} = {\cal F}_{n}(^*)$.  But then ${\cal F}_n(^*) \cup \{v_s v_{s_2}\}= F_n$ implying that $F_n$ is a ${\cal J}_{n}(^*)\/$ configuration of  Type-$1$, a contradiction.

\vspace{1mm}
(ii)  Suppose  ${\cal J}_{n-1}(^*)\/$ is of Type-1.

\vspace{1mm}
 Then there is an edge $e\/$ in ${\cal J}_{n-1}(^*)\/$  such that ${\cal J}_{n-1}(^*) - e \/$  is ${\cal F}_{n-1}(^*)\/$.

 \vspace{1mm} We assert  that $v_{i-1}, v_{i+1}\/$ are both  pendant vertices in $F_n\/$.

\vspace{1mm} To see this   suppose    $v_{i-1}$ is a non-pendant vertex in  ${\cal F}_{n-1}(^*)$. Then it follows from the  maximality of $k\/$ that $v_iv_{i-2} = v_iv_j$ is an edge of $F_n\/$.  But this means that     ${\cal F}_{n-1}(^*) \cup \{v_iv_j\} = {\cal F}_n(^*)$  implying  that    $F_n \/$ is a  $ {\cal J}_n(^*)\/$ configuration of Type-$1\/$ (because  ${\cal F}_{n-1}(^*) \cup \{v_iv_j\} \cup e = {F}_n$), a contradiction.

\vspace{1mm}
Hence let $v_r\/$ and $v_s\/$ be the neighbors of $v_{i-1}\/$ and $v_{i+1}\/$ in $F_n\/$ respectively.

\vspace{1mm}
Suppose first  that either (a) $v_iv_j\/$ crosses neither $v_{i-1}v_r\/$ nor $v_{i+1}v_s\/$ or that (b) $v_iv_j\/$ crosses $v_{i-1}v_r\/$ (or $v_{i+1}v_s\/$) with $|j-r|=1\/$ (or $|j-s| =1\/$). In any of these  cases we have ${\cal F}_{n-1}(^*) \cup \{v_iv_j\} = {\cal F}_n(^*)$ implying that $F_n \/$ is a  $ {\cal J}_n(^*)\/$ configuration of Type-$1\/$, a contradiction.

\vspace{1mm}
 Hence assume without loss of generality  that  $v_iv_j\/$ crosses  $v_{i-1}v_r\/$ with $|j-r| > 1\/$. In this case, let  $G_1$ and $G_2$ be the subgraphs of $K_n - F_n$ induced by the vertices $\{v_{j -1}, v_{j}, \ldots,  v_{i-1}\}$ and $\{v_{i -1}, v_{i},   \ldots, v_{j -1}\}$ respectively.  Clearly  $v_{i-1}$ and $v_i$ are isolated vertex in $G_1\/$ and $G_2\/$ respectively. Hence  $G_i$  admits a triangulation $T_i$ for each $i=1, 2\/$. A  triangulation of $K_n - F_n$ is given by  $T= T_1 \cup T_2$.

\vspace{1mm} This completes the proof.   \qed

\section{Potentially triangulable graphs}     \label{pt}

\vspace{1mm}
We now look at the possibility of packing a graph $F\/$ with  $n\/$ vertices and $n\/$ edges in  the convex complete graph $K_n\/$ so that $K_n - F\/$ admits a triangulation. We shall confine our attention to the case where $F_n\/$ is a $2\/$-regular graph. We begin with the following example.

\vspace{1mm}

\begin{eg}  \label{eg2}
 Suppose  $F\/$ is a  $7\/$-cycle. If $F\/$ is of the form $v_0v_2v_4v_6v_1v_3v_5v_0\/$ or of the form $v_0v_3v_6v_2v_5v_1v_4v_0\/$ in $K_7\/$. Then it is easy to see that $K_7 - F\/$ admits no triangulation. On the other hand, if $F\/$ is of the form $v_0v_2v_6v_4v_1v_3v_5v_0\/$, then $K_7 - F\/$ admits a triangulation whose diagonals are $v_6v_1, v_1v_5, v_5v_2, v_2v_4\/$.
\end{eg}

\begin{define}
Let $K_n\/$ be a convex complete graph with $n\/$ vertices. $F\/$ is said to be   potentially triangulable in $K_n\/$ if there exists a configuration of $F\/$ in $K_n\/$ such that $K_n - F\/$ admits a triangulation.
\end{define}

\begin{result}   \label{fnconnected}
 Suppose $F_n\/$ is an $n\/$-cycle.  Then $F_n\/$ is potentially triangulable in $K_n\/$ if and only if $n \geq 7\/$.
\end{result}

\vspace{2mm}  \noindent
{\bf Proof:}  It is easy to see that $K_n - F_n\/$ admits no triangulation if $n \leq 5\/$.

\vspace{1mm}   Suppose  $K_6 - F_6\/$ admits a triangulation $T\/$. Since $T\/$ has precisely $3\/$ diagonals, there is a vertex $v_i\/$ which is incident with  $2\/$ diagonals of $T\/$. But this is clearly a contradiction since   $K_6 - (F_6 \cup v_0v_1v_2 \cdots v_5v_0)\/$ is $1\/$-regular.

\vspace{1mm} For $n=7\/$, $F_n\/$ is potentially triangulable in $K_n\/$ by Example \ref{eg2}.

\vspace{1mm} For the rest of the proof, we assume that $n \geq 8\/$.

\vspace{1mm}
 {\em Case (1):  $n\/$ is even.}

\vspace{1mm}  Note that $n\/$ can be written as  $n=2^s t\/$ for some positive natural number $s\/$ and some positive odd natural number $t\/$.  Let $\alpha =t+2\/$. Then $gcd \ (t, \alpha) =1\/$.

\vspace{1mm} By relabeling the vertices if necessary, we may assume without loss of generality  that $F_n\/$ takes the form \[ v_0 v_{\alpha} v_{2\alpha} v_{3\alpha} \ldots v_{(n-1)\alpha}v_0.  \]  Here the operation is reduced modulo $n\/$. That is, the edges of $F_n\/$ are of the form $v_i v_j\/$ where $j - i \in S_n = \{t+2, n-t-2\}\/$.

 \vspace{1mm}
 The subgraph $G'\/$ induced by the vertices $v_0, v_2, v_4, \ldots , v_{n-2}\/$  is a convex complete graph $K_{n/2}\/$ (because it contains no edges from $F_n\/$). As such $G'\/$ admits a triangulation $T'\/$ (by Theorem \ref{n-2}). By adding the vertices $v_1, v_3, \ldots, v_{n-1}\/$ to $T'\/$ together with the edges $v_0v_1 v_2v_3 \ldots v_{n-1}v_0\/$, we have a triangulation of $K_n - F_n\/$.

\vspace{1mm}  {\em  Case (2): $n \geq 9\/$ is odd.}

\vspace{1mm}   Let $\alpha = \lfloor n/2 \rfloor\/$ and let $F_n\/$ take the form \[ v_0 v_{\alpha} v_{2\alpha} v_{3\alpha} \ldots v_{(n-1)\alpha}v_0.  \]
Here again,   the operation is reduced modulo $n\/$. That is, the edges of $F_n\/$ are of the form $v_i v_j\/$ where $j - i \in S_n = \{(n-1)/2, (n+1)/2\}\/$.

\vspace{1mm}
Let $\beta = \lceil n/3 \rceil\/$. Consider the subgraphs $G_0, G_1, G_2\/$ induced by the sets of vertices $\{ v_0, v_1, v_2, \ldots, v_{\beta}\}$, $\{v_{\beta}, v_{\beta+1}, v_{\beta+2}, \ldots, v_{2\beta}\}\/$ and $\{v_{2 \beta}, v_{2\beta +1}, \ldots, v_0\}\/$ respectively. Then $G_i\/$, $i=0, 1, 2\/$ is a  convex complete  subgraph of $K_n - F_n\/$ (because it contains no edges from $F_n\/$). Hence $G_i\/$ admits a triangulation $T_i\/$, $i = 0, 1, 2\/$. As such, $T_0 \cup T_1 \cup T_2\/$ is a triangulation of $K_n - F_n\/$.  \qed

\vspace{2mm}

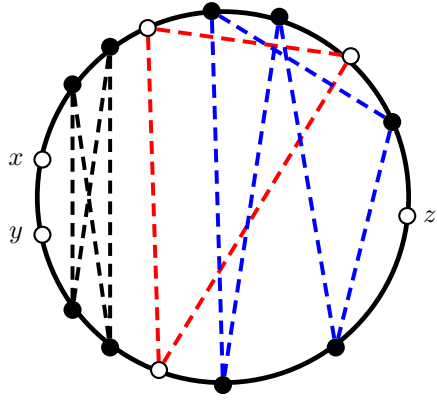
\begin{figure}[htb]
\begin{center}
\psset{unit=0.5cm}
\begin{pspicture}(10,10)

\cnode[linewidth=1.8pt](5,5){2.5cm}{c}

\cnode[fillstyle=solid,fillcolor=white](.2,6){0.12cm}{x}  \cnode[fillstyle=solid,fillcolor=white](.2,4){0.12cm}{y} \cnode[fillstyle=solid,fillcolor=white](9.9,4.5){0.12cm}{z}

\cnode*[linewidth=1.1pt](1,8){0.12cm}{v1}  \cnode*[linewidth=1.1pt](1,2){0.12cm}{u1}
\cnode*[linewidth=1.1pt](2,9){0.12cm}{v2}  \cnode*[linewidth=1.1pt](2,1){0.12cm}{u2}

\ncline[linewidth=1.5pt,linestyle=dashed]{-}{v1}{u1} \ncline[linewidth=1.5pt,linestyle=dashed]{-}{v1}{u2}  \ncline[linewidth=1.5pt,linestyle=dashed]{-}{v2}{u1} \ncline[linewidth=1.5pt,linestyle=dashed]{-}{v2}{u2}

\cnode[fillstyle=solid,fillcolor=white](3,9.5){0.12cm}{v3}    \cnode[fillstyle=solid,fillcolor=white](3.3,.4){0.12cm}{u3}
\cnode[fillstyle=solid,fillcolor=white](8.4,8.75){0.12cm}{v6}


\ncline[linewidth=1.5pt,linecolor=red,linestyle=dashed]{-}{v3}{u3}  \ncline[linewidth=1.5pt,linecolor=red,linestyle=dashed]{-}{v6}{u3}  \ncline[linewidth=1.5pt,linecolor=red,linestyle=dashed]{-}{v3}{v6}

\cnode*[linewidth=2pt](4.7,9.95){0.12cm}{v4}  \cnode*[linewidth=2pt](6.5,9.8){0.12cm}{v5}  \cnode*[linewidth=2pt](9.5,7){0.12cm}{v7}
\cnode*[linewidth=2pt](5,0){0.12cm}{u4}   \cnode*[linewidth=2pt](8,1){0.12cm}{u5}


\ncline[linewidth=1.5pt,linecolor=blue,linestyle=dashed]{-}{v4}{v7}  \ncline[linewidth=1.5pt,linecolor=blue,linestyle=dashed]{-}{v4}{u4}  \ncline[linewidth=1.5pt,linecolor=blue,linestyle=dashed]{-}{v5}{u4}   \ncline[linewidth=1.5pt,linecolor=blue,linestyle=dashed]{-}{v5}{u5}  \ncline[linewidth=1.5pt,linecolor=blue,linestyle=dashed]{-}{v7}{u5}

\rput(-.5,4){\small $y$}   \rput(-.5,6){\small $x$}     \rput(10.5,4.5){\small $z$}

\end{pspicture}
\end{center}

\caption{A drawing of $K_{15} - (C_4 \cup C_3 \cup C_5)\/$.}  \label{xtpolygon}
\end{figure}

\begin{note} \label{r1}
Let $K_n\/$ be a convex complete graph with $n\/$ vertices and let $F\/$ be a union of $k\/$ disjoint cycles  in $K_n\/$ where $|V(F)| =n-3 \geq 4\/$. It is easy to see that one can place these $k\/$ cycles in $K_n\/$ such that (i) no edge of $F\/$ is a boundary edge of $K_n\/$, and (ii) there are three vertices $x, y, z \in V(K_n) - V(F)\/$ such that  $xy\/$ is a boundary edge of $K_n\/$ and  $x, y, z\/$ separate $V(F)\/$ into two sets.   Figure \ref{xtpolygon} illustrates an example where $n=15\/$, $k=3\/$ and $F = C_4 \cup C_3 \cup C_5\/$. Here the cycles are drawn with dashed lines. By Theorem \ref{n-2}, $K_n - F\/$ admits a triangulation (because $|E(F)| = n-3\/$).
\end{note}

\vspace{2mm}
\begin{result} \label{disconnected}
Let $F_n\/$ be a $2\/$-regular graph with $n\/$ vertices. Then $F_n\/$ is potentially triangulable in $K_n\/$ if and only if $F_n \not \in \{ C_3, C_4, C_5, C_6, C_3 \cup C_3, C_3\cup C_4\}\/$.
\end{result}

\vspace{2mm}  \noindent
{\bf Proof:} If $F_n\/$ is connected, the result is true by Theorem \ref{fnconnected}. Hence we assume that $F_n\/$ is a union of disjoint cycles.

\vspace{1mm} Let the vertices of the convex complete graph $K_n\/$ be denoted by $v_0, v_1, \ldots, v_{n-1}\/$.

\vspace{1mm} Since $K_6 - (C_3 \cup C_3  \cup v_0v_1v_2 \cdots v_5v_0)\/$ is $1\/$-regular, it follows that $K_n - F_n\/$ admits no triangulation if $F_n = C_3 \cup C_3\/$.

\vspace{1mm}   Suppose  $F_n =  C_3\cup C_4\/$ and assume that  that $K_n - F_n\/$ admits a triangulation $T\/$.  Since $K_n - F_n\/$ is a $2\/$-regular graph, the graph $D\/$  induced by the the diagonal edges of $T\/$ is a path with $5\/$ vertices. Assume without loss of generality that $D = v_1v_6v_2v_5v_3\/$. It follows that $v_0v_2v_4v_0 \subset F_n\/$. Consequently $E(C_4) \cap E(D) \neq \emptyset\/$, a contradiction.

\vspace{1mm} For the rest of the proof, we assume that $n \geq 8\/$. Let $C\/$ be a cycle in $F_n\/$ and let $|V(C)| = p\/$. Also let $F^* = F_n - C\/$.

\vspace{1mm}
 {\em Case (1):  $p \geq 8\/$.}

\vspace{1mm} (i) Suppose $p\/$ is even.
Use the method in Case (1) of the proof of Theorem \ref{fnconnected} to construct a triangulation $T_p\/$ of $K_p - C\/$. To obtain a triangulation for $K_n - F_n\/$, we do the following.

\vspace{1mm} Assume first that $|V(F^*)| \geq 4\/$.
Insert $n-p\/$ vertices into the edges $v_0v_1, v_1v_2\/$ so that these $n-p\/$ vertices together with $v_0, v_1, v_2\/$ form a  convex complete graph $K_{n-p+3}\/$.  Also, join each of these $n-p\/$ vertices to every vertex in $K_p\/$.    We then place all the disjoint cycles $F^*\/$ in  $K_{n-p+3}\/$ in the same manner as was done in Remark \ref{r1} (with $v_0, v_2, v_1\/$ playing the roles of $x, y, z\/$ respectively) and obtain a triangulation $T^*\/$ for $K_{n-p+3} - F^*\/$.   Consequently, $T_p \cup T^*\/$ is a triangulation for $K_n - F_n\/$.

\vspace{1mm}  Now assume that $F^*\/$ is a $3\/$-cycle.  Insert $3\/$ vertices $u, v, w\/$ into the edges $v_0v_1, v_1v_2, v_2v_3\/$ respectively so that $\{v_0, u, v_1, v, v_2,  w, v_3\}\/$ becomes the vertex set of a convex complete graph $K_7\/$. Also, join $u, v, w\/$ each to every vertex of $K_p\/$. Take $uvwu\/$ to be the $3\/$-cycle $F^*\/$. Then $T_p \cup \{v_0uv_1vv_2wv_3\}\/$ is  a triangulation for $K_n - F_n\/$.

\vspace{1mm} (ii) Suppose $p\/$ is odd.  Use the method in Case (2) of the proof of Theorem \ref{fnconnected} to construct a triangulation $T_p\/$ of $K_p - C\/$. Recall that $T_p = T_0 \cup T_1 \cup T_2\/$ where $T_i\/$ is a triangulation of $G_i\/$, $i =0, 1, 2\/$. To obtain a triangulation for $K_n - F_n\/$, we do the following.

\vspace{1mm} As before we first  assume that $|V(F^*)| \geq 4\/$.
Insert $n-p\/$ vertices into the edges $v_0v_1, v_1v_2\/$ so that these $n-p\/$ vertices together with $v_0, v_1, v_2\/$ form a  convex complete graph $K_{n-p+3}\/$. Also, join each of these $n-p\/$ vertices to every vertex in $K_p\/$.      We then place all the disjoint cycles $F^*\/$ in  $K_{n-p+3}\/$ in the same manner as was done in Remark \ref{r1} (with $v_0, v_2, v_1\/$ playing the roles of $x, y, z\/$ respectively) and obtain a triangulation $T^*\/$ for $K_{n-p+3} - F^*\/$ (since $v_0, v_1, v_2\/$ are isolated vertices in $K_{n-p+3} -F^*\/$). Then  $T_0 \cup T_1 \cup T_2 \cup T^*\/$ is a triangulation for $K_n - F_n\/$.

\vspace{1mm} Now assume that $F^*\/$ is a $3\/$-cycle. Insert $3\/$ vertices $u, v, w\/$ into the edges $v_0v_1, v_1v_2, v_2v_3\/$ respectively so that $\{v_0, u, v_1, v, v_2,  w, v_3\}\/$ becomes the vertex set of a convex complete graph $K_7\/$. Also, join $u, v, w\/$ each to every vertex of $K_p\/$. Take $uvwu\/$ to be the $3\/$-cycle $F^*\/$.
Note that the subgraph induced by  $\{v_0, u, v_1, v, v_2,  w, v_3, \ldots , u_{\beta} \}\/$ (recall that $\beta = \lceil n/3 \rceil\/$) is of the form $K_{\beta + 3} - F^*\/$ which admits a triangulation $T_0^*\/$ (since $v_0, v_1, v_2, \ldots, v_{\beta}\/$ are isolated vertices in $K_{\beta +3} -F^*\/$).     As such $T_0^* \cup T_1 \cup T_2\/$ is a triangulation for $K_n - F_n\/$.

\vspace{2mm}   {\em Case (2): $p = 7\/$.}

\vspace{1mm}
Let $C\/$ be the $7\/$-cycle $v_0v_3v_1v_6v_4v_2v_5v_0\/$. Let $T_p\/$ be the triangulation of $K_7 - F_7\/$ with $4\/$ diagonals $v_0v_2, v_2v_6, v_6v_3, v_3v_5\/$. Adopt the same method of construction as was done in Case (1)(i) to obtain a triangulation $T_p \cup T^*\/$ for $K_n - F_n\/$ if $|V(F^*)| \geq 4\/$ and a triangulation $T_p \cup \{v_0uv_1v v_2wv_3\}\/$ for $K_n - F_n\/$ if $F^*\/$ is a $3\/$-cycle.

\vspace{2mm} {\em Case (3): $p = 6\/$.}

\vspace{1mm} Let $C\/$ be the $6\/$-cycle $v_0v_2v_5v_3v_1v_4v_0\/$.

\vspace{1mm} Assume first that $|V(F^*)| \geq 4\/$.
Insert $n-p\/$ vertices into the edges $v_0v_1, v_0v_5\/$ so that these $n-p\/$ vertices together with $v_0, v_1, v_5\/$ form a  convex complete graph $K_{n-3}\/$.  We then place all the disjoint cycles $F^*\/$ in  this $K_{n-3}\/$ in the same manner as was done in Remark \ref{r1} (with   $v_5, v_1, v_0\/$ playing the roles of $x, y, z\/$ respectively). Also, join the vertices of $F^*\/$ to every vertex in $v_0, v_1, \ldots, v_5\/$ to obtain the convex graph $K_n - F_n\/$. Let  $V(F_1^*)\/$ and $ V(F^*_2)\/$ denote the sets of vertices of $F^*\/$ which are placed on the edges $v_0v_1\/$ and $  v_0v_5\/$ respectively. Consider the subgraphs $G_0\/$ and $G_1\/$  induced by $V(F_1^*) \cup \{v_0, v_1, v_2, v_3\}\/$ and $V(F^*_2)\cup \{v_0, v_5, v_4, v_3\}\/$ respectively. By Theorem  \ref{n-3}, $G_i\/$ admits a triangulation  $T_i\/$, $i=0, 1\/$.   Consequently, $T_0 \cup T_1 \cup \{v_0v_3\}\/$ is a triangulation for $K_n - F_n\/$.

\vspace{1mm} Now suppose $F^*\/$ is a $3\/$-cycle. Here  we take  $F_n = C_6 \cup C_3\/$ where $C_6 = v_0v_3v_8v_5v_2v_6v_0\/$ and $C_3 = v_1v_4v_7v_1\/$. Then a triangulation for  $K_n - F_n\/$ is given by  the set of diagonals $\{v_5v_0,  v_0v_7, v_7v_5,  v_5v_1,  v_1v_3, v_3 v_5\}\/$.

\vspace{2mm} {\em Case (4): $p = 5\/$.}

\vspace{1mm} Let $C\/$ be the $5\/$-cycle $v_0v_2v_4v_1v_3v_0\/$.

\vspace{1mm} Assume first that $|V(F^*)| \geq 4\/$.
Insert $n-5\/$ vertices into the edges $v_0v_1, v_1v_2\/$ so that these $n-5\/$ vertices together with $v_0, v_1, v_2\/$ form a  convex complete graph $K_{n-2}\/$. We then place all the disjoint cycles $F^*\/$  in  this $K_{n-2}\/$ in the same manner as was done in Remark \ref{r1} (with $v_0, v_2, v_1\/$ playing the roles of $x, y, z\/$ respectively).  Also, join all these $n-5\/$ vertices to each of the vertex  in $\{v_0, v_1, \ldots, v_4\}\/$ so that the resulting convex graph is $K_n - F_n\/$. Let $v_0w\/$ be the boundary edge of $K_{n-2}\/$ where $w \in V(F^*)\/$. Let $G^*\/$ denote the subgraph of $K_n - F_n\/$ induced by the vertices $V(F^*) \cup \{v_1, v_2, v_3\}\/$.  Since $v_2\/$ is an isolated vertex in $G^*\/$, $G^*\/$ admits a triangulation $T^*\/$. As such $T^* \cup \{wv_0, v_0v_4,  v_4w, v_3v_4\}\/$ is a triangulation of $K_n - F_n\/$.

\vspace{1mm} Hence assume that $F^*\/$ is a $3\/$-cycle. Here we take  $F_n = C_5 \cup C_3\/$ where $C_5 = v_0v_3v_6v_2v_5v_0\/$ and $C_3 = v_1v_4v_7v_1\/$. Then a triangulation for  $K_n - F_n\/$ is given by  the set of diagonals $\{v_0v_2, v_2v_4, v_4v_0, v_0v_6,  v_4v_6\}\/$.

\vspace{2mm} {\em Case (5): $p = 4\/$.}

\vspace{1mm}  In view of the preceding cases, we may assume that $F_n\/$ is a disjoint union of cycles each of length at most $4\/$.

\vspace{1mm} Suppose $F_n = C_4 \cup C_4\/$. Then  we may take $F_n\/$ to be the two disjoint $4\/$-cycles given by $v_0v_2v_5v_3v_0\/$ and $v_1v_6v_4v_7v_1\/$. Then a triangulation for $K_8 -F_8\/$ is given by the set of diagonals $v_1v_3, v_1v_4, v_1v_5, v_5v_0, v_5v_7\/$.

\vspace{1mm}     Hence we assume that  $|V(F^*)| \geq 6\/$.
Take a convex complete graph $K_{n-1}\/$ and place all the disjoint cycles $F^*\/$ in $K_{n-1}\/$ in the same manner as was done in Remark \ref{r1}. Now insert a new vertex $w\/$ on a boundary edge incident to $z\/$ a\emph{}nd join $w\/$ to all vertices of $K_{n-1}\/$ to obtain a convex complete graph $K_n\/$. Let $C\/$ denote the cycle $xzywx\/$. Let $V(F_1^*)\/$ denote the set of vertices in $F^*\/$ which are placed on the boundary edge $xz\/$, and let  $V(F_2^*)\/$ denote the set of vertices in $F^*\/$ which are placed  on the boundary edge $wy\/$.

\vspace{1mm}
 Let $x_1 \in V(F_1^*)$ and $x_2 \in V(F_2^*)$  such that $x_1x_2$ is not an edge in $F^*\/$.  Partition $K_n - F_n\/$ into three convex graphs $G_1\/$, $G_2\/$ and $G_3\/$  induced by the vertices $\{x_1, \ldots, z\}$, $\{z, \ldots, x_2\}$ and $\{x_2 , \ldots, y, x, \dots, x_1\}$. Then it is easy to see that   $z$  is an  isolated in $G_1$ and $G_2$   and each $x$ and $y$ is an isolated vertex in $G_3$.  Let $T_i\/$ denote a triangulation of $G_i\/$, $i=1, 2, 3\/$. Then a triangulation for $K_n - F_n\/$ is given by $T_1 \cup T_2 \cup T_3\/$.

\vspace{2mm} {\em Case (6): $p = 3\/$.}

\vspace{1mm} In view of the preceding results, we just need to consider the case  $F_n = C_3 \cup C_3 \cup \cdots \cup C_3\/$ and $n \geq 9\/$.

\vspace{1mm} Let $\alpha = n/3\/$ and take $F_n = \{ v_i v_{i+\alpha}v_{i + 2 \alpha}v_i \ | \ i =0, 1, \ldots, \alpha-1\}\/$.

\vspace{1mm} For $ i=0, 1, 2\/$, let $G_i\/$ denote the convex complete subgraph induced by the vertices $v_{i(\alpha-1)}, v_{i(\alpha-1)+1}, \ldots, v_{(i+1)(\alpha-1)}\/$. Also, let $G_3\/$ denote the convex complete subgraph induced by the vertices $v_{3(\alpha-1)}, v_{3\alpha -2}, \ldots, v_0\/$.

\vspace{1mm} Each $G_i\/$, $i=0, 1, 2, 3\/$ admits a triangulation $T_i\/$ if $n \geq 12\/$. Hence
$T_0 \cup T_1 \cup T_2 \cup T_3 \cup \{v_0v_{\alpha -1}v_{2(\alpha-1)}v_{3(\alpha-1)}v_0, v_0v_{2(\alpha-1)}\}\/$ yields a triangulation for $K_n -F_n\/$ if $n \geq 12\/$. If $n=9\/$, a triangulation for $K_9 -F_9\/$ is given by the diagonals $\{v_0v_2, v_2v_4, v_4v_0,  v_4v_6, v_6v_8, v_8v_4\}\/$.

\vspace{1mm} This completes the proof.   \qed

\section{Regular graphs}  \label{regular}

In view of the results in the preceding section, it is natural to ask which regular graph is potentially  triangulable in $K_n\/$.

\begin{problem}
Let $r \geq 3\/$ be a natural number and let $G\/$ be an $r\/$-regular graph with $n\/$ vertices. It is true that there is a natural number $n_0(r)\/$ such that when $n \geq n_0(r)\/$, then $G\/$ is potentially triangulable in the convex complete graph $K_n\/$?
\end{problem}

\vspace{1mm} We believe that the above problem is true. However we do not have a complete  answer for this even when restricted to the case $r=3\/$.  Nevertheless  we offer the following special case of a $3\/$-regular graph which is well-known in the literature.

\vspace{1mm}
Suppose $n\/$ and $k\/$ are two integers such that $1 \leq k \leq n-1\/$ and $n \geq 5\/$. The {\em generalized Petersen graph\/} $P(n,k)\/$ is defined to have vertex-set
  $\{a_i, b_i : i = 0, 1, \ldots, n-1\}\/$ and edge-set $E_1 \cup E_2 \cup E_3\/$ where   $E_1 = \{ a_ia_{i+1} : i = 0, 1, \ldots, n-1\}\/$,
$E_2 = \{ b_ib_{i+k} : i = 0, 1, \ldots, n-1\}\/$   and $E_3 = \{ a_ib_i : i = 0, 1, \ldots, n-1\}\/$     with subscripts reduced modulo $n \/$. Edges in $E_3\/$ are called the spokes of $P(n,k)\/$.

\begin{propo}
Suppose $ 1 \leq k < n/2\/$. Then the generalized Petersen graph $P(n, k)\/$ is potentially triangulable in the convex complete graph $K_{2n}\/$  where  $n \geq 5\/$.
\end{propo}

\vspace{1mm}  \noindent
{\bf Proof:}  Let the vertices of $K_{2n}\/$ be denoted  $v_1, v_2, v_3, \ldots, v_{2n}\/$.  We shall pack $P(n,k)\/$ on $K_{2n}\/$ so that $K_{2n} - P(n,k)\/$ admits a triangulation.

\vspace{1mm} {\em Case (1): \ $k=1\/$}

\vspace{1mm} In this case, $P(n,1)\/$ consists of two $n\/$-cycles $C = a_0a_1a_2 \cdots c_{n-1}a_0\/$ and $C' = b_0b_1b_2 \cdots \linebreak b_{n-1}b_0\/$ together with the edges $a_ib_i\/$, $i=0, 1, 2, \ldots, n-1\/$.

\vspace{1mm}  If $n\/$ is even, we place  $C\/$ on $K_{2n}\/$ so that $C\/$ takes the form

\vspace{1mm}  \hspace{8mm} $v_1v_{3n/2}v_2v_{(3n-2)/2}v_3v_{(3n-4)/2} \cdots v_{(n-2)/2}v_{n+2}v_{n/2}v_{n+1}v_1\/$

\vspace{1mm} \noindent  and that $C'\/$ takes the form

\vspace{1mm}  \hspace{8mm}    $v_{(3n+2)/2} v_n v_{(3n+4)/2} v_{n-1}v_{(3n+6)/2} v_{n-2} \cdots v_{(n+4)/2} v_{2n} v_{(n+2)/2} v_{(3n+2)/2} \/$

\vspace{1mm} \noindent
with  spokes given by \ $v_{(3n+2i)/2}v_i\/$, \ $i =1, 2, \ldots, n/2\/$ \ and \ $v_{(3n-2i)/2} v_{n-i}\/$, \ $i = 0, 1, 2, \ldots, n/2-1\/$.

\vspace{1mm}  If $n\/$ is odd, we place  $C\/$ on $K_{2n}\/$ so that $C\/$ takes the form

\vspace{1mm}  \hspace{8mm} $v_1 v_{(3n-1)/2} v_2 v_{(3n-3)/2} v_3 v_{(3n-5)/2} \cdots v_{n+2}v_{(n-1)/2} v_{n+1}v_{(n+1)/2} v_1\/$

\vspace{1mm} \noindent  and that $C'\/$ takes the form

\vspace{1mm}  \hspace{8mm} $v_{(3n+1)/2} v_n v_{(3n+3)/2} v_{n-1} v_{(3n+5)/2} v_{n-2} \cdots v_{2n-1}v_{(n+3)/2} v_{2n} v_{(3n+1)/2}\/$

\vspace{1mm} \noindent
with  spokes given by \  $v_{(3n+2i-1)/2} v_i\/$, \ $i=1, 2, \ldots, (n+1)/2\/$, \ and \ $v_{(3n-2i-1)/2} v_{n-i}\/$, \ $i=0, 1, 2, \ldots, (n-3)/2\/$.

\vspace{1mm} In both cases consider the subgraph graphs $G_1\/$ and $G_2\/$   induced by the sets of vertices $\{v_1, v_2, \ldots, v_n\}\/$ and $\{v_{n+1}, v_{n+2}, \ldots, v_{2n-1},  v_{2n}\}\/$ respectively. Since  the vertex $v_2\/$ in $G_1\/$ (likewise the vertex $v_{n+1}\/$ in $G_2\/$) is not incident with any edges from $G_1\/$ (respectively  $G_2\/$) since  $n \geq 5\/$, $G_i\/$ admits a triangulation $T_i\/$, $i=1, 2\/$. A  triangulation of $K_{2n} - P(n,1)\/$  is given by $T_1 \cup T_2 \cup \{ v_nv_{2n}, v_1v_{2n}, v_nv_{n+1}\}\/$.

\vspace{2mm} {\em Case (2): \ $1 < k < n/2\/$}

\vspace{1mm} Place $E_1\/$ and $E_2\/$ on $K_{2n}\/$ so that  $E_1\/$ takes the form $v_2v_4v_6 \ldots v_{2n-2}v_{2n}v_2\/$,  $E_2\/$ takes the form $\{ v_iv_{i+2k} \ : \ i =1, 3,  \ldots, 2n-1\}  \/$ and $E_3\/$ takes the form $\{v_{2i}v_{2i+3} \ : \ i= 1, 2, \ldots, n\}\/$. Here the operations on the subscripts are reduced modulo $2n\/$.

\vspace{1mm} Consider the subgraph $G\/$ induced by the vertices $\{v_4, v_5, v_6, \ldots, v_{2n-3}\} \cup \{ v_0, v_1, v_{2n-1}\}\/$. Since the vertex $v_0\/$ is not adjacent to every vertex in $G\/$, $G\/$ admits a triangulation $T_G\/$. Then $T_G \cup \{v_1v_2v_3v_4 \} \cup \{v_1v_3\} \cup \{ v_{2n-3}v_{2n-2}v_{2n-1}\}\/$ is a triangulation for $K_{2n} - P(n,k)\/$.

\vspace{1mm} This completes the proof.  \qed

\vspace{1mm}

\section{An application} \label{apply}

The skewness of a graph $G\/$, denoted $sk(G)\/$, is the minimum number of edges to be deleted from $G\/$ such the resulting graph is planar.   In \cite{kai1:refer}, Kainen proves that all graphs $G\/$ with $sk(G) \leq 2\/$ are $5\/$-colorable (and the result is best possible). Also, in the same paper he proves that if $sk(G) \leq 5\/$, then $G\/$ is $6\/$-colorable.  In the following year,  in \cite{kai2:refer}, it was shown that if $sk(G) < { k \choose 2}\/$, then $G\/$ is $(2+k)\/$-colorable for $k \geq 3\/$. This result was  generalized to other orientable surfaces in \cite{pah:refer}.

\vspace{1mm}
For more details  about the notion of skewness of a graph and for some recent results, the reader may consult the papers \cite{cl2:refer}, \cite{cs:refer} and \cite{cim:refer}.

 \vspace{1mm}


\begin{define} \label{convexlemma}
 The geometric skewness of a geometric graph  $G\/$, denoted $sk_g(G)\/$ is the minimum number of edges to be removed from $G\/$ so that the resulting graph can be redrawn as a geometric planar graph.  The convex skewness of a convex graph $G\/$, denoted $sk_c(G)\/$ is the minimum number of edges to be removed from $G\/$ so that the resulting graph is a convex plane graph.
\end{define}

\begin{result}
For any geometric graph $G\/$, $sk_g(G) = sk(G)\/$.
\end{result}

\vspace{2mm} \noindent
{\bf Proof:} The proof follows immediately from F\'{a}ry's theorem (\cite{far:refer}) which states that any simple planar graph can be drawn on the plane without crossings so that its edges are straight line segments. \qed


\vspace{1mm}
\begin{propo} \label{triangle}
Let $F\/$ be a subgraph of a convex complete graph $K_n\/$. Suppose $K_n-F\/$ admits a triangulation. Then $sk_c(K_n-F) = {n-2\choose 2} - |E(F)|\/$.
\end{propo}

\vspace{2mm} \noindent
{\bf Proof:}   It is known that the number of edges in a triangulation $T\/$ of a convex $n\/$-gon  is $2n-3\/$ (with  $n-3\/$ of them being non-boundary edges). If any new straight line segment is added to the triangulation, it will intersect with some non-boundary edge of $T\/$. Hence, if  $K_n - F\/$ admits a triangulation, then  we have $sk_c(K_n- F) = |E(K_n)| - |E(F)| - (2n-3)\/$ which yields  $sk_c(K_n - F) = {n-2\choose 2} - |E(F)|\/$.  \qed

\end{document}